\documentclass{eptcs}
\providecommand{\event}{ADG 2021} 
\usepackage{graphicx}
\usepackage{amssymb} 
\usepackage{mathtools} 

\def\orcidID#1{\unskip$^{\mbox{\href{https://orcid.org/#1}{\scriptsize{[#1]}} }}$}

\begin{document}

\title{Parametric Root Finding \\for Supporting Proving and Discovering \\ Geometric Inequalities in GeoGebra}
\author{
Zolt\'an Kov\'acs\orcidID{0000-0003-2512-5793}
\institute{
The Private University College of Education of the Diocese of Linz\\
Linz, Austria}
\email{zoltan@geogebra.org}
\and 
R\'obert Vajda\orcidID{0000-0002-2439-6949}
\institute{Bolyai Institute, University of Szeged  \\
Szeged, Hungary}
\email{vajda@math.u-szeged.hu}
}

\def\titlerunning{Parametric Root Finding for$\ldots$ Proving and Discovering Geometric Inequalities in GeoGebra}
\def\authorrunning{Kov\'acs \& Vajda}

\maketitle              

\begin{abstract}
We introduced the package/subsystem \emph{GeoGebra Discovery} to GeoGebra 
which supports the automated proving or discovering of elementary geometry inequalities.     
In this case study, for inequality exploration problems related to isosceles and right angle triangle subclasses,
we demonstrate how our general real quantifier elimination (RQE) approach could be replaced by a parametric root finding (PRF) algorithm.
The general RQE requires the full cell decomposition of a high dimensional space, while the new method can avoid this expensive computation and can lead to practical speedups. To obtain a solution for a 1D-exploration problem, 
we compute a Gr\"obner basis for the discriminant variety of the 1-dimensional parametric system and solve finitely many nonlinear real (NRA) satisfiability (SAT)  problems.
We illustrate the needed computations by examples. Since Gr\"obner basis algorithms are available in Giac (the underlying free computer algebra system in GeoGebra) and freely available efficient NRA-SAT solvers 
(SMT-RAT, Tarski, Z3, etc.) can be linked to GeoGebra, we hope that the method could be easily added to the existing reasoning tool set for educational purposes.  
\end{abstract}


\section{Introduction}

As we reported in our earlier paper \cite{KVR} and recent works \cite{BKVT}, 
the dynamic geometry system GeoGebra \cite{GGBA1, GGBA2} supports an automated reasoning toolset (ART).
In particular, a GeoGebra user may try to prove or explore a relation between geometric quantities defined by a standard (planar) Euclidean
construction.  That is, for instance, we may wish to prove that in each non-degenerate triangle, the ratio of the sum of the medians and the perimeter of a triangle cannot exceed 1.
To phrase it differently, we may want to explore general \emph{elementary geometric inequalities} related to Euclidean plane geometry constructions.    
Recently, based on general real quantifier elimination (RQE)\cite{Coll1}, the \emph{realgeom} \cite{VK18} tool supports this automated exploration in the background.

However, due to high complexity of the general RQE/CAD algorithms \cite{DH} and the big  number of variables in the input RQE formula, some of the well known and elementary results \cite{BOT} are inaccessible with our approach and implemented tools.

In this case study we want to replace the solution method based on general purpose RQE with a new approach where we may avoid the full cylindrical decomposition of a high dimensional space at least for a particular problem class. We hope that with the new method speedups can be obtained.


\section{Parametric Real Root Finding (PRF)}
We will consider only two simple subclasses of elementary geometry inequality exploration challenges:
\begin{itemize}
\item problems related to isosceles triangles (IT) and
\item for right angle triangles (RT).
\end{itemize}

The reasons are as follows. The proving/disproving of IT/RT-conjectures of inequality types via algebraic methods, after the algebraic formulation leads to a real nonlinear satisfiability (NRA-SAT) problem. That is, to an existentially closed formula which validity should be decided. 

In contrast, for the IT/RT exploration problems we are concerned with in this paper, the typical associated first order input formula contains \emph{one free variable}, and the semialgebraic system is (generically) \emph{one-dimensional}: the quantities which we want to compare, does not have a fixed ratio in a triangle, but it can vary from triangle to triangle from the investigated class IT (or RT). 
Still, in some sense, the translated algebraic problems are very close to NRA-SAT problems: for each fixed parameter $m=m_0$, the system has finitely many (maybe zero) real solutions. Thus we can avoid a generic real quantifier elimination process for determining the possible range of the parameter where a real solution exists via constructing  a full CAD of the $r$-dimensional space, where $r$ is defined by the number of the variables in  the input formula. 

Instead, by knowing the  Discriminant Variety (DV), which characterizes the ``critical/wrong" points $W$, and which could be determined by (the hopefully cheaper) Gr\"obner basis computations, we      
can be sure that the system behaves well (it has constantly many solutions) in the open cells of a 1D-decomposition $\mathbb{R}\backslash W$ determined by the DV-polynomials.
Therefore, sampling the open cells and adding the finitely many zero-dimensional cells in $W$ to the sample, we may reduce the exploration problem to finitely many NRA-SAT (in fact, to nonparametric real root counting (RRC)) problems.
      
We call this method the \emph{parametric root finding} method (PRF). For the the details of the parametric real root counting and for the role of the (minimal) discriminant variety the reader is referred to \cite{Lazard, Liang}.

\section{Examples}
\subsection{Example 1: A problem  from IT}
In GeoGebra, for each construction and to the related exploration problem, the translated semi-algebraic problem is based on the planar coordinates 
of the geometric points involved in the construction of the objects (vertices, midpoints of triangles, intersection points of lines, etc.)  
However, in this first example we intentionally avoid yet the introduction of the coordinatization, to reduce the number of variables in the input formula and making the related computational steps as simple as possible. We note that the number of indeterminates in the automatically derived input formulae based on coordinatization is much higher, typically between 4-10. However, the number of parameters in these problems,
no matter if they were derived by coordinatization or they are coordinate-free versions,
 is constant, namely, 1.

\begin{figure}[h]
\begin{center}
\includegraphics[scale=0.6]{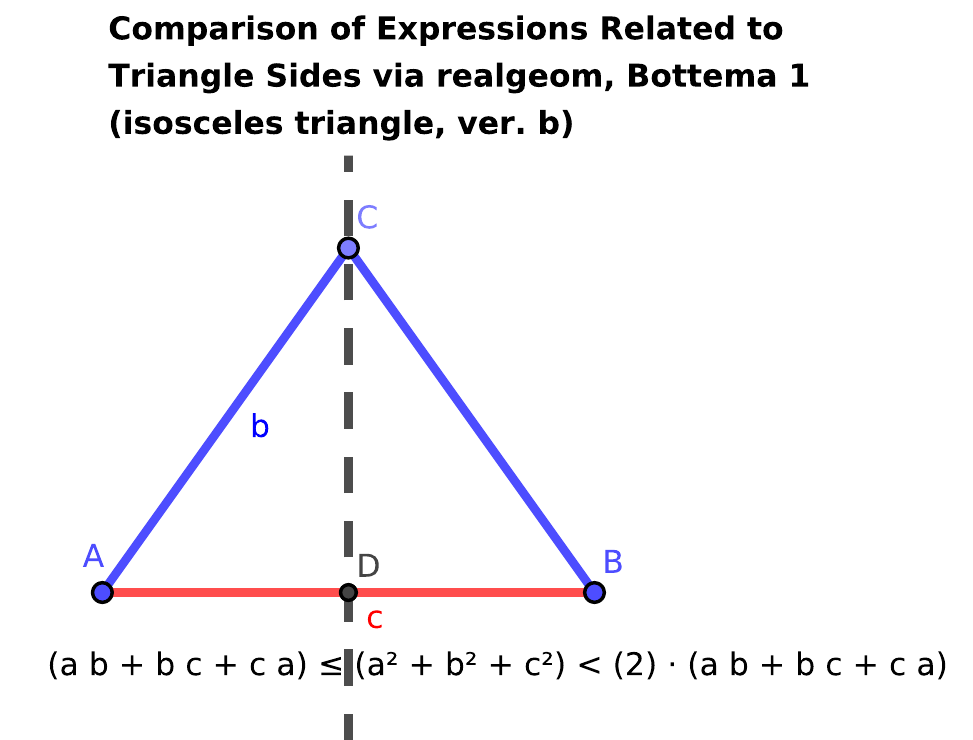}
\includegraphics[width=5cm]{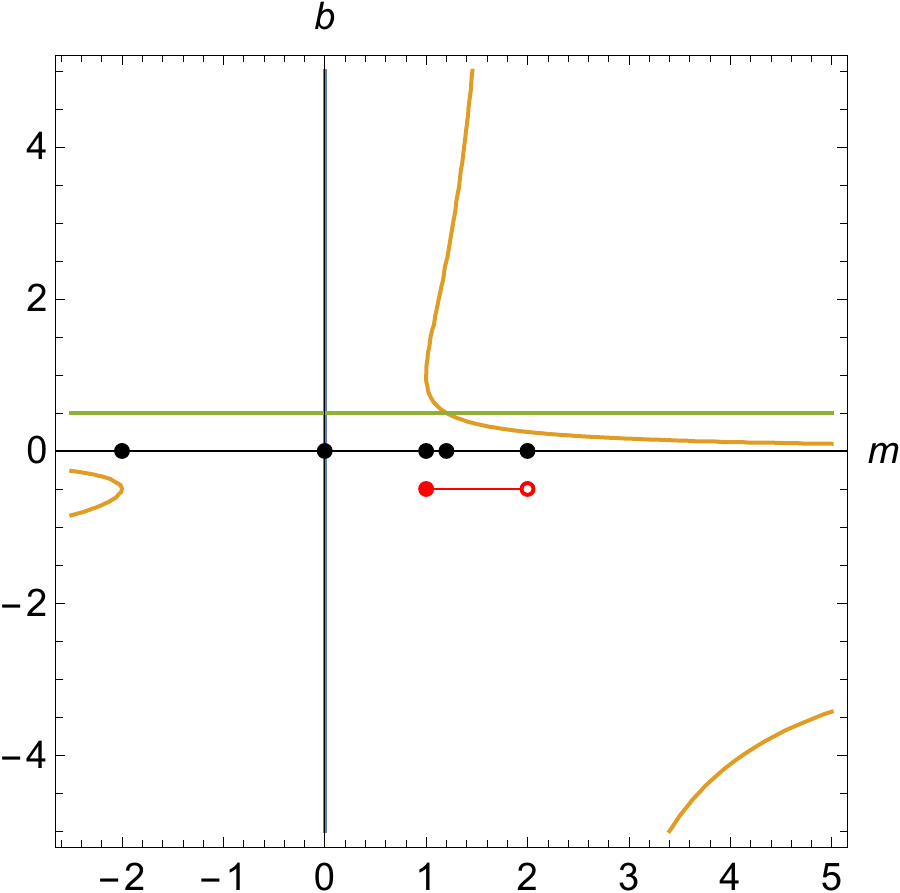}
\caption{An exploration problem for a non-degenerate isosceles triangle with GeoGebra Discovery and the $(m,b)$-space}\label{figure:1} 
\end{center}
\end{figure}

Assume that an isosceles triangle with vertices $A,B,C$ and side lengths $a=b,b,c$ is given (see Figure \ref{figure:1}, left). 
Without loss of generality we can assume that $c=1$. 
We want to investigate (explore) the range of the ratio $m$ of 
\begin{equation}
\frac{AB^2+BC^2+CA^2}{AB\cdot BC+AB\cdot CA+ BC\cdot CA}=\frac{1+2b^2}{b^2+2b}.
\end{equation} 

The related standard input formula for RQE would look like 
\begin{equation}\label{inputformula1} 
\underset{b}{\exists}\;{m>0 \wedge 2b-1>0 \wedge (1+2b^2)=m\, (b^2+2b)},
\end{equation} 
and the quantifier free output is $1\le m<2$.
We show now that the same output formula can be computed with a different/new method.
As a first step we compute the univariate polynomials for the discriminant variety $DV$
\begin{equation}
O_{{\rm crit}}=(2+m)(m-1),\quad O_{{\rm in}}=m (-6 + 5 m), \quad O_{{\rm inf}}=2-m. 
\end{equation}
To obtain the univariate polynomials in $O_{\rm crit}, 
O_{\rm in}, O_{\rm inf}$, we follow \cite[Theorem 2]{Moroz6}, but for an update see \cite{Moroz11}.
Let $f= (1+2b^2)-m (b^2+2b), g_1=2b-1, g_2=m$. 

For $O_{crit}$ we need the determinant of the partial Jacobian $J(f)$ (w.r.t.~$b$):
\begin{equation}
J(f)=\frac{\partial (1+2b^2)-m(b^2+2b)}{\partial b}=b(4-2m)-2m.
\end{equation}     
Then we compute a Gr\"obner basis of the elimination ideal of $\left<f, J(f), t \cdot g_1\cdot g_2-1\right> \cap\, \mathbb{Q}[m]$ and get $m^2+m-2$.

In a similar way, for $O_{\rm in}$, that is,  for inequalities, we compute the Gr\"obner basis of the ideal $(\left< f, g_1\cdot g_2-u, u t-1\right> \cap\, \mathbb{Q}[m,u])_{u=0}$ and we get as factors $m$ and $5m-6$.

The remaining univariate polynomial(s) for $O_{\rm inf}$ are also computed via Gr\"obner bases, but besides the elimination and specialization, we need homogenizations as well. We skip the details here.

The real positive roots of the polynomials in $DV$ in increasing order are $d_2=1,d_4=6/5,d_6=2$. That is, $W=\{-2,0,1,6/5,2\}$.
We choose sample points from the open cells $d_1=(0,1),d_3=(1,6/5),d_5=(6/5,2), d_7=(2,\infty)$ and add them to the positive roots:
\begin{equation}
\{1/2,1,11/10,6/5,3/2,2,3\}.
\end{equation}
With NRA-SAT (RRC) we obtain that for the cells $d_2, d_3, d_4, d_5$  there is at least one real solution of the above semialgebraic system when $m$ is replaced by the sample in the cell.
Therefore we conclude that $m=1 \vee 1<m<6/5 \vee m=6/5 \vee 6/5<m<2\equiv 1\le m<2$ is the solution for the exploration problem. See Figure \ref{figure:1}, right.

\subsection{Example 2: A problem from RT}

\begin{figure}[h]
\begin{center}
\includegraphics[scale=0.4]{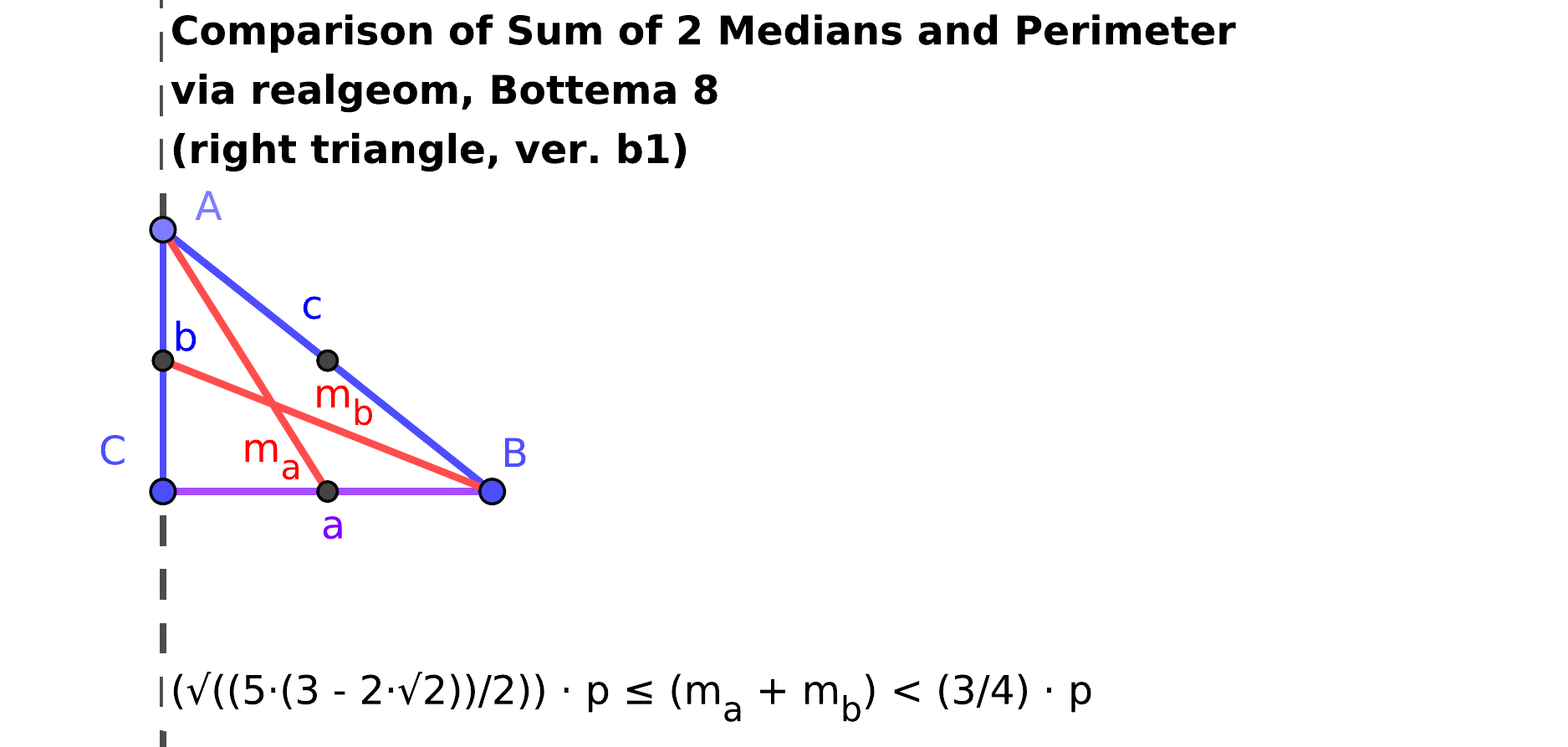}
\caption{Example 2: An exploration problem for a non-degenerate right triangle with GeoGebra Discovery}\label{figure:2} 
\end{center}
\end{figure}

Our second worked example considers a non-degenerate right angle triangle with vertices $A, B, C$ and hypotenuse $c=AB$, where the exploration task is to investigate the ratio of sum of the medians $m_a$ and $m_b$ and the perimeter $p=a+b+c$ (see Figure \ref{figure:2}). 

In the first order input formula now we still have the parameter $m$, but instead of $a,b,c,m_a,m_b$ we will see the (bound) variables 
$v_{13}, v_{15}, v_{16}, v_{17}$. These variables were generated in GeoGebra Discovery in a mechanical way: the input first order formula was automatically derived, based on the coordinates of the triangle vertices and midpoints of the triangle sides. Some indexed $v_j$'s, where indexing starts from 1, are missing, because they were (linearly) eliminated in a preprocessing step or they were set to special values without loss of generality:    
\begin{align}\label{inputformula2} 
\underset{v_{13},v_{15},v_{16},v_{17}}{\exists}\;{m>0 \wedge v_{13}>0 \wedge v_{15}>0 \wedge v_{16}>0 \wedge v_{17}>0\, \wedge}\\
 v_{13} + v_{16} - m (1 + v_{15} + v_{17}) = 0 \wedge
 15 + 4 v_{13}^2 - 16 v_{16}^2 = 0\, \wedge \nonumber\\
 3 - 4 v_{16}^2 + v_{17}^2 = 0 \wedge 4 + v_{15}^2 - 4 v_{16}^2 = 0.\nonumber
\end{align} 
For the reader's convenience, we note that the variables $v_{15}$ and $v_{17}$ correspond to the triangle side lengths $b$ and $c$, and the 
variables $v_{13}$ and $v_{16}$ to the triangle medians $m_a$ and $m_b$.  
The quantifier free output is $\sqrt{\frac{5}{2} (3-2\sqrt2)}\le m<3/4$.
Now this time the univariate polynomials in DV look like 
\begin{align}
O_{{\rm crit}}&=(25 - 60 m^2 + 4 m^4)\\ 
&(32805 - 523422 m^2 + 388800 m^4 + 
   1377792 m^6 + 737280 m^8 + 131072 m^{10}),\nonumber \\
O_{{\rm in}}&=m (-3 + 4 m) (-1 + 4 m) (1 + 4 m) (3 + 4 m) (-3 + 4 m^2)\\
& (3 + 
   4 m^2) (50625 - 324000 m^2 + 633600 m^4 + 368640 m^6 + 65536 m^8),\nonumber\\
O_{{\rm inf}}&=(-3 + 4 m) (-1 + 4 m) (1 + 4 m) (3 + 4 m), 
\end{align}
and therefore we have to work also with higher order algebraic numbers in our sample:
\begin{align}
d_2=1/4, d_4=r^{(10,3)}, d_6=r^{(10,4)}, d_8=\sqrt{\frac{5}{2} (3-2\sqrt2)},\\
 d_{10}=3/4, d_{12}=\sqrt{3}/2, d_{14}=\sqrt{\frac{5}{2} (3+2\sqrt2)}, \nonumber
\end{align}
where $r^{(10,3)}$ and $r^{(10,4)}$ refer to the third and fourth positive roots of the degree 10 polynomial in $O_{{\rm crit}}$. 
This time only two consecutive SAT problems (with samples from $d_8$ and $d_9$) evaluate to True and the final formulae obtained is
\begin{equation}
m=\sqrt{\frac{5}{2} (3-2\sqrt2)} \vee \sqrt{\frac{5}{2} (3-2\sqrt2)}<m<3/4.  
\end{equation}

\section{Discussion/Conclusion}\label{Disc}

Variants of the general parametric real root finding and real root counting algorithms are implemented and available in Maple \cite{ChenMaza, Liang}, but to our best knowledge, they are not available directly in free computer algebra systems. Also our prototype implementation for the  1D explorations is in Mathematica \cite{MMA1}. Therefore, in our future work we will make statistics of the required computational times for the Gr\"obner basis computation and for solving the NRA-SAT problems related to the harder exploration problems in IT/RT.

As a preliminary example, in Table \ref{table:1} we give the number of cells reported by QEPCAD B, version 1.72 \cite{qepcad}. The cells were constructed when solving a 1-dimensional, six variable exploration problem related to the ratio of the inradius and circumradius ($m=R/r$) in right triangles, via CAD-based RQE.
\begin{table}
\label{table:1}
\begin{center}
 \begin{tabular}{c c c c c c c} 
 Level & 1 & 2 & 3 & 4 & 5 & 6\\
 \hline
 Cells & \;175\; & 379\; & \;89008\; & \;40644\; & \;541\; & \;529\;
\end{tabular}
\end{center}
\caption{Number of cells in QEPCAD B, version 1.72}
\end{table}
Here the generated output formula is $m\ge\sqrt{2}+1$ and the DV-polynomials generated by the PRF method are relatively simple:
\begin{equation}
 DV=\{m, 2m-1,2m+1, m^2-2m-1, m^2+2m-1\}. 
 \end{equation}
We see that the the construction of the more than $10^5$ cells may be replaced by 7 pieces of 1D cells in this example.
 
 If we gain practical speedups, then we intend to adapt and re-implement the PRF method and solve the computational subproblems with Giac and the freely available NRA-SAT solvers. Thus the emerging tool may be used in a broader educational context.

Note that the question about the possible ratio
\begin{equation}
m=\frac{AB^2+BC^2+CA^2}{AB\cdot BC+AB\cdot CA+ BC\cdot CA},
\end{equation}
which is discussed in Example 1, could have been also asked about a nondegenerate general triangle as an exploration problem. 
We observed that most of the exploration problems for general triangles lead to semialgebraic problems, which have still infinitely many solutions for a fixed $m=m_0$ (after setting one of the triangle side lengths to 1). 
However, considering another variable in the input formula, say, one of the triangle sides $b$, as an additonal parameter to $m$, that is, if we fix $b=b_0$ and $m=m_0$, then the number of solutions of the semialgebraic system is again finite.
 
This opens up the road to an application of the same idea and methods for solving some more difficult geometric exploration problems without a full CAD. However, after investigating the open cells of the 2D space, we would have to cover here also all the infinitely many parameter pairs where the bivariate polynomials in the DV vanish. This can be done by adding the bivariate polynomials one by one to the original semialgebraic system and performing new ``real root parametric counting" for the new systems. An additonal problem may occur, because some resulting, again 1-dimensional semialgebraic system may have solutions with multiplicity $>1$ for almost all parameter values; in this case the DV-based method cannot be applied directly  (cf.~\cite[Theorem 2]{Moroz6}). Still, we think the above general idea may be elaborated in the future for handling geometric exploration problems for arbitrary triangles as well.    


\section{Acknowledgements}
The first author was partially supported by the grant PID2020-113192GB-I00 (Mathematical Visualization: Foundations, Algorithms and Applications) from the Spanish MICINN.
The second author was supported by the EU-funded Hungarian grant EFOP-3.6.1-16-2016-00008.

\providecommand{\urlalt}[2]{\href{#1}{#2}}
\providecommand{\doi}[1]{doi:\urlalt{http://dx.doi.org/#1}{#1}}

\end{document}